\documentclass[a4paper,11pt]{amsart}

\usepackage[T1]{fontenc}

\usepackage{amssymb}
\usepackage{amsmath}

\usepackage[utf8]{inputenc}

\usepackage[english,francais]{babel}

\newtheorem{theorem}{Theorem}[section]
\newtheorem{lemma}[theorem]{Lemma}
\newtheorem{e-proposition}[theorem]{Proposition}

\newtheorem{e-definition}[theorem]{Definition\rm}
\newtheorem{remark}{\it Remark\/}

\newcommand{\R}{\mathbf{R}}
\newcommand{\abs}[1]{\vert #1  \vert}
\DeclareMathOperator{\Div}{div}

\def\og{\leavevmode\raise.3ex\hbox{$\scriptscriptstyle\langle\!\langle$~}}
\def\fg{\leavevmode\raise.3ex\hbox{~$\!\scriptscriptstyle\,\rangle\!\rangle$}}

\begin{document}

\title[Pathological solutions]{Pathological solutions to elliptic problems in divergence form with continuous coefficients}
\address{Rutgers University,
Department of Mathematics,
110 Frelinghuysen Road,
Piscataway NJ 08854-8019, United States of America}
\author{Tianling Jin}
\email{kingbull@math.rutgers.edu}
\author{Vladimir Maz'ya}
\address{University of Liverpool, Department of Mathematical Sciences, 
Liverpool L69 3BX, United Kingdom}
\address{Link\"oping University, Department of Mathematics,
SE-581 83 Link\"oping, Sweden}
\email{vlmaz@mai.liu.se}
\author{Jean Van Schaftingen}
\address{Universit\'e catholique de Louvain, D\'epartement de Math\'ematique, Chemin du Cyclotron 2, 1348 Louvain-la-Neuve, Belgium}
\email{Jean.VanSchaftingen@uclouvain.be}

%


\maketitle




\begin{abstract}
\selectlanguage{english}%
We construct a function $u \in W^{1,1}_{\mathrm{loc}} (B(0,1))$ which is a solution to $\Div (A \nabla u)=0$ in the sense of distributions, where $A$ is continuous and $u \not \in W^{1,p}_{\mathrm{loc}} (B(0,1))$ for $p > 1$. We also give a function $u \in W^{1,1}_{\mathrm{loc}} (B(0,1))$ such that $u \in W^{1,p}_{\mathrm{loc}}(B(0,1))$ for every $p < \infty$, $u$ satisfies $\Div (A \nabla u)=0$ with $A$ continuous but $u \not \in W^{1, \infty}_{\mathrm{loc}}(B(0,1))$.

\end{abstract}

\section{Introduction}
\label{}

Consider the equation
\begin{align}
\label{eqn}
  -\Div A \nabla u &= 0 & & & \text{in $\Omega$},
\end{align}
for $\Omega \subset \R^n$. If  $A  : \Omega \to \R^{n \times n}$ is bounded, measurable and elliptic, i.e., there exists $\lambda, \Lambda \in \R_*$ such that for every $x \in \Omega$, $A(x)$ is a symmetric matrix, and 
\[
  \abs{\xi}^2 \le (A(x)\xi)\cdot \xi \le \Lambda \abs{\xi}^2,
\] then one can define a weak solution $u \in W^{1,1}_{\mathrm{loc}}(\Omega)$ by requiring that for every $\varphi \in C^1_c(\Omega)$,
\[
 \int_{\Omega} (A \nabla u)\cdot \nabla \varphi = 0.
\]

We are interested in the regularity properties of $u$. 
A fundamental result of E. De Giorgi \cite{DG} states that if $u \in W^{1,2}_{\mathrm{loc}}(\Omega)$, then $u$ is locally H\"older continuous. In particular, $u$ is then locally bounded. In the same direction, N. G. Meyers \cite{M} also proved that $u \in W^{1,p}_{\mathrm{loc}}(\Omega)$ for some $p > 2$. 

J. Serrin \cite{S} showed that the assumption $u \in W^{1,2}_{\mathrm{loc}}(\Omega)$ is essential in E. De Giorgi's result by constructing for every $p \in(1, 2)$ a function $u \in W^{1,p}_{\mathrm{loc}}(\Omega)$ that solves such an elliptic equation but which is not locally bounded. In these counterexamples $A$ is not continuous. J. Serrin \cite{S} conjectured that if $A$ was H\"older continuous, then any weak solution $u \in W^{1,1}_{\mathrm{loc}}(\Omega)$ is in $W^{1,2}_{\mathrm{loc}}(\Omega)$, and one can then apply E.~De~Giorgi's theory.

This conjectured was confirmed for $u \in W^{1,p}(\Omega)$ by R. A. Hager and J. Ross \cite{HR} and for $u \in W^{1,1}(\Omega)$ by H. Brezis \cite{Ancona,BLincei}. The proof of Brezis extends to the case where the modulus of continuity of $A$ 
\begin{equation}
\label{omegaA}
 \omega_A(t)=\sup_{\substack{x, y \in \Omega\\ \abs{x-y} \le t}} \abs{A(x)-A(y)},
\end{equation}
satisfies the Dini condition
\begin{equation}
\label{Dini}
 \int_0^1 \frac{\omega_A(s)}{s}\,ds < \infty.
\end{equation}

In the case where $A$ is merely continuous, H. Brezis obtained the following result

\begin{theorem}[H. Brezis \cite{Ancona,BLincei}]
Assume that $A \in C(\Omega; \R^{n \times n})$ is elliptic. If $u \in W^{1,p}_{\mathrm{loc}}(\Omega)$ solves~\eqref{eqn}, then for every $q \in [p, +\infty)$, one has $u \in W^{1,q}_{\mathrm{loc}}(\Omega)$.
\end{theorem}

H. Brezis asked two questions about the cases $p=1$ and $q=\infty$ in the previous theorem. 
We answer both questions, with a negative answer. First we have

\begin{e-proposition}
\label{propositionW11}
There exists $u \in W^{1,1}_{\mathrm{loc}}(B(0,1))$ and an elliptic $A \in C(B(0,1); \R^{n \times n})$ such that $u$ solves~\eqref{eqn}, but $u \not \in W^{1,p}_{\mathrm{loc}}(B(0,1))$ for every $p > 1$.
\end{e-proposition}

As a byproduct, we obtain

\begin{e-proposition}
\label{propositionNonunique}
There exists $A \in C(B(0,1); \R^{n \times n})$ such that the problem
\begin{equation}
\label{eqBVP}
  \left\{
\begin{aligned}
   -\Div (A \nabla u) &= 0 &&\text{in $B(0,1)$},\\
   u&=0 &&\text{on $\partial B(0,1)$}.
\end{aligned}
  \right.
\end{equation}
has a nontrivial solution.
\end{e-proposition}

Our construction in Proposition~\ref{propositionW11} allows in fact to show that the counterexamples can be improved

\begin{e-proposition}
\label{propositionW1log}
There exists $u \in W^{1,1}_{\mathrm{loc}}(B(0,1))$ and an elliptic $A \in C(B(0,1); \R^{n \times n})$ such that $u$ solves~\eqref{eqn}, $Du \in (L \log L)_{\mathrm{loc}} (B(0,1))$ but $u \not \in W^{1,p}_{\mathrm{loc}}(B(0,1))$ for every $p > 1$.
\end{e-proposition}

In particular, in this case, $Du$ belongs locally to the Hardy space $\mathcal{H}^1$ (see \cite{Stein1969}).

%

Concerning the possibility of Lipschitz estimates, we have

\begin{e-proposition}
\label{propositionW1infty}
There exists $u \in W^{1,1}_{\mathrm{loc}}(B(0,1))$ and an elliptic $A \in C(B(0,1); \R^{n \times n})$ such that $u$ solves~\eqref{eqn}, $Du \in W^{1,p}_{\mathrm{loc}}(B(0,1))$ for every $p > 1$, $Du \in \mathrm{BMO}_{\mathrm{loc}}(B(0,1))$ but $u \not \in W^{1, \infty}_{\mathrm{loc}}((B(0,1))$.
\end{e-proposition}

This shows that  $Du \in L^p(B(0,1))$ does not imply $Du \in L^\infty(B(0,1/2))$, one can wonder whether it implies that $Du \in \mathrm{BMO}(B(0,1/2))$. The answer is still negative

\begin{e-proposition}
\label{propositionBMO}
There exists $u \in W^{1,1}_{\mathrm{loc}}(B(0,1))$ and an elliptic $A \in C(B(0,1); \R^{n \times n})$ such that $u$ solves~\eqref{eqn}, $u \in W^{1,p}_{\mathrm{loc}}(B(0,1))$ for every $p \in (1, \infty)$  but $Du \not \in \mathrm{BMO}_{\mathrm{loc}}(B(0,1))$.
\end{e-proposition}

The construction of the counterexamples are made by explicit formulas, inspired by the construction of J. Serrin \cite{S}. They can also be obtained from asymptotic formulas of V.~Kozlov and V.~Maz'ya \cite{KM2003a,KM2003b}.

\section{The pathological solutions}

Our counterexamples rely on the following computation

\begin{lemma}
\label{lemma}
Let $v \in C^2((0,R))$ and $\alpha \in C^1((0,R))$. Define $A(x)=(a_{ij}(x))_{\substack{1 \le i \le n \\ 1 \le j \le n}}$ by
\[
 a_{ij}(x)=\delta_{ij}+\alpha(\abs{x})\Bigl(\delta_{ij}-\frac{x_i x_j}{\abs{x}^2}\Bigr).
\]
Then for every $x \in B(0,R) \setminus \{0\}$, 
\begin{equation}
\label{identity}
  \Div \big(A(x) \nabla (x_1 v(\abs{x}))\big)
=x_1 \Bigl(v''(\abs{x})+\frac{n+1}{\abs{x}} {v'(\abs{x})}-\frac{n-1}{\abs{x}^2}\alpha(\abs{x})v(\abs{x})\Bigr).
\end{equation}
\end{lemma}


\begin{remark}
If $P$ is a homogeneous harmonic polynomial of degree $k$, the formula generalizes to
\begin{multline}
  \Div \big(A(x) \nabla (P(x)v(\abs{x}))\big)\\
=P(x) \Bigl(v''(\abs{x})+\frac{n+2k-1}{\abs{x}}v'(\abs{x})-\frac{k(n+k-2)}{\abs{x}^2}\alpha(\abs{x})v(\abs{x})\Bigr).
\end{multline}
\end{remark}

\begin{proof}[Proof of Proposition \ref{propositionW11}]
Choose $\beta > 1$, and define for some $r_0 > 1$, for $r \in (0, 1)$, 
\begin{equation}
\label{vW11}
 v(r)=\frac{1}{r^{n}(\log \frac{r_0}{r})^\beta}.
\end{equation}
One takes then
\begin{equation}
\label{alphaW11}
  \alpha(r)=\frac{r^2 v''(r)+(n+1)rv'(r)}{(n-1) v(r)}
=\frac{-\beta n}{(n-1)\log \frac{r_0}{r}}+\frac{\beta(\beta+1)}{(n-1)\bigl( \log \frac{r_0}{r} \bigr)^2}.
\end{equation}
One can take $r_0$ large enough so that $\alpha \ge -\frac{1}{2}$ on $(0,1)$; the coefficient matrix $A$ is then uniformly elliptic. Define now $u(x)=x_1 v(\abs{x})$. One checks that $u \in W^{1,1}(B(0,1))$ and that $u$ is a weak solution of~\eqref{eqn}. Indeed, it is a classical solution on $B(0,1) \setminus \{0\}$ by the previous lemma. Taking, $\varphi \in C^1_c(B(0,1))$ and $\rho \in (0,1)$, and integrating by parts we obtain
\[
\begin{split}
 \int_{B(0,1)\setminus B(0,\rho)} \nabla \varphi \cdot (A \nabla u)&=
-\int_{\partial B(0,\rho)} \varphi \nabla u \cdot (A \frac{x}{\rho})\\
&=-\int_{\partial B(0,\rho)} \varphi \nabla u \cdot \frac{x}{\rho}\\
&=-\int_{\partial B(0,\rho)}\varphi x_1 \Big( \frac{v(\rho)}{\rho}+v'(\rho)\Big)  \\
&=-\int_{\partial B(0,\rho)} (\varphi(x)-\varphi(0)) x_1 \Big( \frac{v(\rho)}{\rho}+v'(\rho)\Big).
\end{split}
\]
Since $\varphi \in C^1_c(B(0,1))$, one has
\[
\Bigl\lvert \int_{B(0,1)\setminus B(0,\rho)} \nabla \varphi \cdot (A \nabla v)\Bigr \rvert \le 
C \rho^{n}(\abs{v(\rho)}+\rho\abs{v'(\rho)}),
\]
since the right-hand side goes to $0$ as $\rho \to 0$, $u$ is a weak solution.
\end{proof}

\begin{remark}
The examples constructed in the case of merely measurable coefficients by J. Serrin \cite{S} to show that a solution $u \in W^{1,p}_{\mathrm{loc}}(\Omega)$ need not be in $W^{1,2}_{\mathrm{loc}}(\Omega)$ and by N. G. Meyers \cite{M} to show that for every $p > 2$, that a solution in $W^{1,2}_{\mathrm{loc}}(\Omega)$ need not be in $W^{1,p}_{\mathrm{loc}}(\Omega)$ can be recovered with the same construction, by taking $v(r)=r^\alpha$. The ellipticity condition requires $\alpha < n-1$ or $\alpha >  1$. This covers all the cases when $n=2$; a descent argument finishes the construction in higher dimension.
\end{remark}

\begin{proof}[Proof of Proposition \ref{propositionW1log}]
One checks that the counterexample constructed in the proof of Proposition~\ref{propositionW1log} satisfies $Du \in L \log L (B(0,1))$ when $\beta > 2$. 
\end{proof}

Similar examples can be obtained following the results of V.~Kozlov and V.~Maz'ya \cite{KM2003b}. By (4) therein, if $A \in C(B(0,1); \R^{n \times n})$, $A(Rx)=R A(x) R$ where $R$ is the reflection with respect to the $x_1$ variable and $A$ satisfies some regularity assumptions, then the equation $-\Div (A \nabla u)=0$ has a solution that is odd with respect to the $x_1$ variable and that behaves like 
\[
 \frac{x_1}{\abs{x}^n}\exp \Bigl(\int_{B(0,1) \setminus B(0, \abs{x})} \mathcal{R}(y)\,dy\Bigr)
\]
around $0$, where $\mathcal{R}$ is defined following \cite[(3)]{KM2003b}\footnote{The reader should correct the misprint in \cite[(3)]{KM2003b} and read $\abs{S^{n-1}_+}$ instead of $\abs{S^{n-1}}$ .}
\begin{equation}
\label{eqR}
  \mathcal{R}(x)= \frac{(e_1 \cdot (A(x)-A(0)) e_1) (x \cdot A(0)^{-1} x)-n ( e_1 \cdot (A(x)-A(0)) A(0)^{-1} x) (e_1 \cdot x)}{\abs{\partial B(0,1)} \abs{\det A(0)}^\frac{1}{2} ( x \cdot A(0)^{-1} x)^{\frac{n}{2}+1}}.
\end{equation}
Taking $A$ as in Lemma~\ref{lemma} with $\lim_{r \to 0} \alpha(r)=0$, one has $\mathcal{R}(x)=\alpha(\abs{x})(\abs{x}^2-x_1{}^2)/(\abs{\partial B(0,1)} \abs{x}^{n+2})$. Therefore, there is a solution that behaves like 
\[
\frac{x_1}{\abs{x}^n} \exp \bigl( \frac{n-1}{n} \int_{\abs{x}}^1 \alpha(r)\,\frac{dr}{r}\bigr).
\]
In particular, if one takes $\alpha(r)=-\beta n/((n-1)\log \frac{r_0}{r})$, one obtains a solution that behaves like $\frac{x_1}{\abs{x}^n} (\log \frac{r_0}{r})^{-\beta}$.
One could also take $a_{ij}(x)=\delta_{ij} +\kappa(\abs{x}) (\delta_{ij}- n \delta_{i1}\delta_{j1}\frac{x_1{}^2}{\abs{x}^2})$ and continue the computations with now $\mathcal{R}(x)=\kappa(\abs{x}) (\abs{x}^2-n x_1{}^2)^2/(\abs{\partial B(0,1)} \abs{x}^{n+2})$.

\begin{proof}[Proof of Proposition \ref{propositionNonunique}]
Let $u$ be given by the proof of Proposition \ref{propositionW11}. Notice that $u$ is smooth on $\partial B(0,1)$. Since $A$ is bounded and elliptic, the problem 
\[
  \left\{
\begin{aligned}
   -\Div (A \nabla w) &= 0 &&\text{in $B(0,1)$},\\
   w&=u &&\text{on $\partial B(0,1)$}.
\end{aligned}
  \right.
\]
has a unique solution in $w \in W^{1,2}(B(0,1))$. Since $u \not \in W^{1,2}(B(0,1))$, $u \ne w$. Hence, $u-w \in W^{1,1}(B(0,1))$ is a nontrivial solution of~\eqref{eqBVP}.
\end{proof}

\begin{proof}[Proof of Proposition \ref{propositionW1infty}]
Take for $r \in (0, 1)$, 
\begin{equation}
\label{vW1infty}
 v(r)=\log \frac{r_0}{r}
\end{equation}
and
\begin{equation}
\label{alphaW1infty}
 \alpha(r)=\frac{1-(n+1)}{(n-1)\log \frac{r_0}{r}}=\frac{-n}{(n-1)\log\frac{r_0}{r}},
\end{equation}
where $r_0$ is chosen so that $\alpha(r) > -\frac{1}{2}$ on $(0, 1)$.
Defining $u(x)=x_1 v(\abs{x})$, one checks that $Du \in W^{1,p}_{\mathrm{loc}}(B(0,1))$, $Du \in \mathrm{BMO}(B(0,1))$,  $u \not \in W^{1,\infty}(B(0,1))$ and that $u$ solves~\eqref{eqn} in the sense of distributions.
\end{proof}

As for the previous singular pathological solutions, similar examples can be obtained from results of V.~Kozlov and V.~Maz'ya for solutions \cite{KM2003a}. By (4) therein if $A \in C(B(0,1); \R^{n \times n})$, $A(Rx)=R A(x) R$ where $R$ is the reflection with respect to the $x_1$ variable and $A$ satisfies some regularity assumptions, then the equation $-\Div (A \nabla u)=0$ has a solution in $W^{1,2}(B(0,1))$ that is odd with respect to to the $x_1$ variable and that behaves like 
\[
  x_1 \exp \Bigl(-\int_{B(0,1) \setminus B(0, \abs{x})} \mathcal{R}(y)\,dy\Bigr)
\]
around $0$, where $\mathcal{R}$ is given by~\eqref{eqR}.
Taking $A$ as in Lemma~\ref{lemma} with $\alpha(r)=\frac{-n}{n-1}(\log \frac{r_0}{r})^{-1}$ one recovers the counterexample presented above.

\begin{proof}[Proof of Proposition \ref{propositionBMO}]
Define for  $r \in (0, 1)$, 
\[
 v(r)=\big(\log \frac{r_0}{r}\big)^2.
\]
and
\[
 \alpha(r)=\frac{-2n}{(n-1)\log\frac{r_0}{r}}+\frac{2}{(n-1)(\log \frac{r_0}{r})^2}.
\]
Defining $u(x)=x_1 v(\abs{x})$, one checks that $u \in W^{1,p}(B(0,1))$ for every $p > 1$ and that $u$ solves~\eqref{eqn} in the sense of distributions.
One checks that for every $c > 0$, $\exp (c \abs{Du}) \not \in L^1(B(0, \frac{1}{2}))$; hence by the John--Nirenberg embedding theorem \cite{JN} (see also e.g. \cite[Chapter 4, \S 1.3]{Stein}), $Du \not \in \mathrm{BMO}(B(0,\frac{1}{2}))$.
\end{proof}

\section*{Acknowledgments}

The authors thank H. Brezis for bringing their attention on the problem. J.V.S. aknowledges the hospitality of the Mathematics Department of Rutgers University.

\end{document}